\documentclass[12pt,oneside, english]{article}
\usepackage[english]{babel}
\usepackage{amsmath, amssymb, amsthm, amscd, color}
\usepackage{cancel}
\usepackage{cite}

\usepackage{enumerate}
\newtheorem{theorem}{Theorem}[section]
\newtheorem{definition}[theorem]{Definition}

\newtheorem{remark}[theorem]{Remark}
\newtheorem{lemma}[theorem]{Lemma}
\newtheorem{corollary}[theorem]{Corollary}
\newtheorem{proposition}[theorem]{Proposition}

\def\fq{\mathbb{F}_{q}}
\def\fqs{\mathbb{F}_{q^2}}
\def\fqn{\mathbb{F}_{q^n}}

\def\cF{\mathcal{F}}

\def\cC{\mathcal{C}}
\def\cF{\mathcal{F}}
\def\cP{\mathcal{P}}
\def\cG{\mathcal{G}}
\def\cH{\mathcal{H}}
\def\cL{\mathcal{L}}
\def\cS{\mathcal{S}}
\def\cN{\mathcal{N}}

\title{Locally recoverable codes from automorphism groups of function fields of genus $g \geq 1$}
\author{Daniele Bartoli\thanks{Daniele Bartoli is with the Dipartimento di Matematica e Informatica, Universit\`a degli Studi di Perugia, Via Vanvitelli 1 - 06123 Perugia - Italy, {\em email:} daniele.bartoli@unipg.it} \, Maria Montanucci \thanks{Maria Montanucci is with the Department of Applied Mathematics and Computer Science, Technical University of Denmark, Kongens Lyngby 2800, Denmark,  {\em email: } marimo@dtu.dk}, 
and Luciane Quoos\thanks{Luciane Quoos is with the Instituto de Matem\'atica, Universidade Federal do Rio de Janeiro, Rio de Janeiro 21941-909 - Brazil, {\em email: } luciane@im.ufrj.br}
}

\date{}

\begin{document}

\maketitle

\begin{abstract}
A Locally Recoverable Code is a code such that the value of any single coordinate of a codeword can be recovered from the values of a small subset of other coordinates. When we have $\delta$ non overlapping subsets of cardinality $r_i$ that can be used to recover the missing coordinate we say that a linear code $\cC$ with length $n$, dimension $k$, minimum distance $d$ has $(r_1,\ldots, r_\delta)$-locality  and denote by
$[n, k, d; r_1, r_2,\dots, r_\delta].$ In this paper we provide a new upper bound for the minimum distance of these codes. Working with a finite number of subgroups of cardinality $r_i+1$ of the automorphism group a function field $\cF| \fq$ of genus $g \geq 1$ we propose a construction of $[n, k, d; r_1, r_2,\dots, r_\delta]$-codes and apply the results to some well known families of function fields.
\end{abstract}   

{\bf Keywords:} Codes for distributed storage, Locally Recoverable Codes, Hamming distance, covering maps, maximal curves. \\
\indent{\bf MSC 2010 Codes:} 14G50, 11T71, 94B65. \\

\section{Introduction}
The study of Locally Recoverable codes (LRC for short) was motivated by the use of coding theory techniques applied to distributed and cloud storage systems. Local recovery techniques enable to repair lost encoded data by a local procedure, that is by making use of small amount of data instead of all information contained in a codeword. Formally, a LRC code of length $n$ is a code that produces an $n$-symbol codeword from $k$ information symbols and, for any symbol of the codeword,  there exist at most $r$ other symbols such that the value of the symbol can be recovered from them. This value $r$ is the called the locality of the code. For example, a code of length $2k$ in which each coordinate is repeated twice, is an LRC code with locality $r = 1$. Generally the locality parameter satisfies $1\leq  r \leq k$ since the entire codeword can be found by accessing $k$ symbols other than the erased symbol. 

In distributed storage systems, erasure codes with locality $r$ are preferred because a coordinate can be locally repaired by accessing at most $r$ other coordinates. However, the local repair may not be performed when some of the $r$ coordinates are also erased. To overcome this problem, we can work with  $\delta $ non overlapping local repair sets of size no more than $r_i$ for a coordinate. The formal definition of a Locally Recoverable Code is given in Definition \ref{locality}. We denote  a linear  code $\cC$ with length $n$, dimension $k$, minimum distance $d$, and $(r_1,\ldots, r_\delta)$-locality  by 
$[n, k, d; r_1, r_2,\dots, r_\delta].$

In recent years, the study of locally reparable codes has attracted a lot of attention. Most of the results concern bounds on the minimum distance \cite{GHSY2012, TB2014bound, WZ2014, KNF2019}  and construction of LRC codes \cite{JKZ2019, HMM2018, JKZ, LLX2019, LXC2019, MT2018, MTT2019, TB2014, WZ2014, Micheli}. 

The goal of this paper is twofold. First we present a new upper bound on the minimum distance of a $[n, k, d; r_1, r_2,\dots, r_\delta]$-code. Second we present two constructions of Locally Recoverable Codes over function fields $\cF| \fq$ of genus $g \geq 1$ using subgroups of their automorphism group. We apply these constructions on some well known families of curves with many rational points.

Let  $\cC$ be an $[n, k, d; r_1, \ldots, r_\delta]$-code. If $\delta=1$ it was proved in 2012 that the minimum distance of the code is upper bounded by
\begin{equation}\label{rec1}
d \leq n-k-\left \lceil \frac{k}{r} \right\rceil +2;
\end{equation}
see \cite{GHSY2012}.
This bound coincides with the classical Singleton bound when $r=k$. A code that achieves equality in \eqref{rec1} is called an optimal LRC code. In \cite{TB2014, LXC2019,LLX2019}  the authors  constructed optimal LRC codes using particular types of polynomials, cyclic codes, and elliptic curves respectively. 

For general $\delta \geq 1$ and $r=r_1= \cdots = r_\delta$ the bound (\ref{rec1}) was generalized in 2014 to 
\begin{equation}\label{rect}
d \leq n-\sum_{i=0}^t\left \lfloor \frac{k-1}{r^i} \right\rfloor;
\end{equation}
see \cite{TB2014}.

Since in our paper we deal with general $[n, k, d; r_1, \ldots, r_\delta]$-codes, in order to check the quality of the constructed codes, we first generalize \eqref{rec1} to
\begin{equation}\label{MainBound}
d \leq n-k- \left\lceil\frac{(k-1)\delta+1}{1+\sum_{i=1}^{\delta} r_i}\right\rceil +2;
\end{equation}
see Theorem \ref{Th:Bound}. 
In analogy with the definition of the singleton defect for linear codes, a notion of relative defect for $[n, k, d; r_1, \ldots, r_\delta]$-code is introduced using the bound in \eqref{MainBound} as
\begin{equation}\label{Eq:RSD}
\delta(\mathcal{C}) =\frac{1}{n}\left(n-k-d+2-\left\lceil\frac{(k-1)\delta+1}{1+\sum_{i=1}^{\delta}r_i}\right\rceil\right).
\end{equation}
Clearly, the smaller the relative defect of a code, the  better the code.  

In this work we also provide two general constructions of LRC codes with locality $(r_1, \ldots, r_\delta)$ over function fields $\cF| \fq$ of genus $g \geq 1$ using  a finite number of subgroups of cardinality $r_i+1$ of the automorphism group $Aut(\cF| \fq)$. 
We distinguish the cases of trivial intersection (Theorem  \ref{Th:GeneralConstruction}) and nontrivial intersection (Theorem \ref{Th:GeneralConstructionNon-Trivial}) and we provide a generalization of the construction proposed in \cite{JKZ2019} for the rational function field.

It is well known that maximal function fields and function fields with many rational points  provide algebraic geometric codes with ``good" parameters; see e.g. \cite{BMZ2018, BQZ2018, CT2016,  HZ2016}. We build up several families of LRC codes applying  Theorems  \ref{Th:GeneralConstruction}  and \ref{Th:GeneralConstructionNon-Trivial} to families of maximal function fields, such as the Hermitian \cite{GSX2000}, the Giulietti-Korchm\'aros \cite{GK2009}, the generalized Hermitian \cite{KKO2001, HKK2006}, and the Norm-Trace  \cite{MTT2008} ones.  In particular, we obtain families of locally repairable codes with good relative parameters, since the relative defect $\delta(\mathcal{C})$ as in \eqref{Eq:RSD}  tends to zero when $q$ goes to infinity; see Propositions \ref{herm1}, \ref{herm2}, \ref{GK}, \ref{genherm}, and \ref{normtrace}. Also, we compare the parameters of some of these LRC codes with the ones constructed in \cite{HMM2018} obtained from fiber products of algebraic curves; see Remark \ref{Remark:Comparison}.

The paper is organized as follows. In Section 2, we present some preliminaries on locally repairable codes, function fields, automorphism groups, and algebraic geometry codes. In Section 3
 a new upper bound for the minimum distance $d$ of a $[n, k, d; r_1, r_2,\dots, r_\delta]$-code depending on the locality $(r_1, r_2,\dots, r_\delta)$ is provided; see Theorem \ref{Th:Bound}. This bound is used to define the relative defect $\delta$ of such code, a parameter that will be used to measure how good a code is.  In Sections 3 and 4 we present  two types of constructions of locally repairable codes via subgroups of the automorphism group of a function field. In Theorem \ref{Th:GeneralConstruction} we deal with subgroups with trivial intersection, and in Theorem \ref{Th:GeneralConstructionNon-Trivial} the case of non-trivial intersection is investigated. In the same sections we apply the constructions to some well known function fields with many rational places. 

\section{Preliminaries}

Consider an $[n, k, d]$-code $\mathcal{C}$ over $\fq^n$ and generator matrix $G=(g_1,\ldots, g_n)$. For any subset $I\subset \{1,\ldots, n\}$ let $\langle g_j\,:\, j \in I\rangle$ denote the $\fq$-subspace generated by the set $\{g_j\, :\, j \in I\}$. 

The code $\cC \subseteq \fq^n$ is said to have locality $(r_1, r_2, \dots , r_\delta)$ if the value of every coordinate of a given codeword can be retrieved by accessing any of the $\delta$ disjoint recovering sets of cardinalities $r_1, r_2, \dots, r_\delta$,  i.e., one error in the codeword can be corrected in a local way if no more than $\delta-1$ erasures occur. A formal definition of a locally repairable code is given as follows.

\begin{definition}\cite{BTS2017}\label{locality}
The $i$-th coordinate of a given codeword, $1 \leq i \leq n$, of an $[n, k, d]$ linear code $\mathcal{C}$ whose generator matrix is $G=(g_1,\ldots, g_n)$ is said to have $(r_1,\ldots, r_\delta)$-locality if there exist pairwise disjoint repair sets $R_1^{(i)},\ldots,R_{\delta}^{(i)} \in \{1,\ldots, n\}\setminus \{i\}$  such that for each $1\leq j\leq \delta$
\begin{enumerate}[i)]
\item $\#R_j^{(i)} = r_i$;
\item  $g_i \in \langle g_\ell\ \mid \ell \in R_j^{(i)}\rangle$.
\end{enumerate}
\end{definition}
Note that $r_i = 1$ implies repetition and we only consider codes with $ \delta\geq 1$ and $r_i\geq 2$. Additionally, we always assume $r_i < k$. We denote  a linear  code $\cC$ with length $n$, dimension $k$, minimum distance $d$, and with  $(r_1,\ldots, r_\delta)$-locality  by 
$$[n, k, d; r_1, r_2,\dots, r_\delta].$$

A set $I\subset \{1,\ldots, n\}$ is called an {\it information set} for an $[n, k, d; r_1, r_2,\dots, r_\delta]$-code if $\#I =rank(\langle g_j\mid j \in I\rangle)=k$

%\footnote{Daniele and Maria: the definition is this one and it is used in the proof of the theorem on the bound}.
%In terms of generators matrices we can use the following equivalent definition of locality.

We now fix some notation and collect some results on function fields and algebraic geometry codes; for more on this subject we refer to the book \cite{S2009}. Given a  function field $\cF | \fq$ and a function $z \in \cF$  let $(z)_\infty$ be its pole divisor. For a divisor $D$ in $\cF$ the $\fq$-vector space
$$\cL(D)=\{ z \in \cF \, : \, (z)_\infty \geq - D\} \cup \{0\}$$
is the Riemann-Roch space associated to $D$ and its dimension is denoted by $\ell(D)$. 
Let 
$$Aut(\cF | \fq)=\{\sigma: \cF \mapsto \cF \,:\, \sigma \text{ is an automorphism of } \cF \text{ and } \sigma(a)=a, \, \forall a \in \fq \}$$ 
be the automorphism group of the function field $\cF | \fq$. For a subgroup $\cH$ of $Aut(\cF | \fq)$ we let $\cF^\cH$ denote the fixed field of $\cH$. The group $Aut(\cF | \fq)$ acts on set of rational places of the function field in a natural way and, for any  rational place $P \in \cF$ and any  subgroup $\cH$ of $Aut(\cF | \fq)$, we let $P^\cH=\{\sigma(P) \,:\, \sigma \in \cH \}$ stand for the orbit of $P$ under the action of $\cH$. The orbit is either short or long provided  $P^\cH <\cH$ or $P^\cH=\cH$ respectively. 

Fixing $\cH$ and $P\in \cF^\cH$ an $\fq$- rational place, the ramification of  $P$ in the Galois extension $\cF|\cF^\cH$ is determined by the action of $\cH$ on the orbit  
$P_1^\cH$ where $P_1$ is any fixed $\fq$- rational place over $P$. In particular, $P$ is totally ramified if and only if $P_1^\cH = \{ P_1\}$ , and  $P$ is completely split
 if and only if $\# P_1^\cH= \# \cH$ and every place $P_1$ above $P$ is $\fq$-rational.

Let $G$ be a divisor on $\cF$ and $P_1, P_2, \dots, P_n$ be pairwise distinct rational places on $\cF$, with $P_i\notin supp(G)$ for all $i$.  Define $D=\sum_{i=1}^n P_i$. The linear algebraic geometry code $\cC_\cL(D, G)$ is defined as the image of the evaluation function
$$ev\, : \, \cL(G) \rightarrow \fqn, \, f \mapsto (f(P_1), f(P_2),\ldots, f(P_n)).$$
The $\cC_\cL(D, G)$ code has length $n$ and the classical bound on the minimum distance $d$ is 
\begin{equation}\label{goppabound}
d \geq n-\deg(G).
\end{equation}

\section{A bound on the minimum distance}

Let $\cC$ be an $[n, k, d; r_1, \dots, r_\delta]$-code. Locally recoverable (LRC-)codes attaining equality in \eqref{rec1} are called {\it optimal}. Constructions of optimal LRC-codes
on curves of genus one can be found in \cite{LLX2019}, and with minimum distance $3$ and $4$ where obtained using cyclic codes in \cite{LXC2019}. 
%For general $\delta \geq 1$ and $r=r_1= \cdots = r_\delta$ the bound on (\ref{rec1}) on the minimum distance was generalized in \cite{TB2014} to 
%\begin{equation}\label{rect}
%d \leq n-\sum_{i=0}^\delta\left \lfloor \frac{k-1}{r^i} \right\rfloor.
%\end{equation}

In \cite[Theorem 1]{WZ2014} another upper bound on  the minimum distance was derived for the case $\delta \geq 1$ and $r=r_1= \cdots = r_\delta$
\begin{equation}\label{rec3}
d \leq n-k+2-\left \lceil \frac{(k-1)\delta+1}{(r-1)\delta+1} \right\rceil.
\end{equation}
In what follows we generalize the ideas in \cite{WZ2014}  in order to provide a new bound on the minimum distance of a general $[n, k, d; r_1, \dots, r_\delta]$-code. First we recall a well known lemma about the minimum distance of a code.

\begin{lemma}\cite{MacWilliams}\label{Lemma:MinimumDistance}
The minimum distance of any $[n, k, d]$-code whose generator matrix is $G=(g_1,\ldots, g_n)$ satisfies 
$$d = n - max\{\#N \  : \ N\subset\{1,\ldots ,n\}, \ rank(\langle g_j\mid j \in N \rangle) < k\}.$$
\end{lemma}

%%%%%%%%%%%%%%%%%%%%%%%%%%%%%%%%%%%%%%%%%%%%%%%%%%%
%%%%%%%%%%%%%%%%	MAIN BOUND		%%%%%%%%%%%%%%%%%%%%%%
%%%%%%%%%%%%%%%%%%%%%%%%%%%%%%%%%%%%%%%%%%%%%%%%%%%
\begin{theorem}\label{Th:Bound}
For an $[n, k, d;r_1,\ldots,r_{\delta}]_q$ linear code, with $2\leq r_1\leq r_2\leq \cdots \leq r_{\delta}<k$,
\begin{equation}\label{Eq:MainBound}
d \leq n-k+2- \left\lceil\frac{(k-1)\delta+1}{1+\sum_{i=1}^{\delta}r_i}\right\rceil.
\end{equation}
\end{theorem}
\proof
Let $I\subset \{1,\ldots, n\}$  be an information set such that each coordinate in $I$ has $(r_1,\ldots, r_{\delta})$-locality. For any $i \in I$ and $j=1, \dots ,\delta$ let $R_j^{(i)}$ be the repair set of cardinality $\#R_j^{(i)} = r_i$, and $1\leq a \leq \delta$ denote by $N_a^{(i)}$ the set $\{i\} \cup R_1^{(i)}\cup \cdots \cup R_a^{(i)}$. So, 
\begin{equation}\label{rankN}
rank(\langle g_j\ : \ j \in N_a^{(i)}\rangle)\leq 1+\sum_{i=1}^{a}r_i.
\end{equation}
Since $I$ is an information set there exists $\{i_1,\ldots, i_{\ell}\}\subset I$ such that 

$$rank(\langle g_j\ : \ j \in \overline{N} \rangle)= k-1,$$
where 
$$\overline{N}=N_{\delta}^{(i_1)}\cup \cdots \cup N_{\delta}^{(i_{\ell-1})}\cup R_1^{(i_{\ell})}\cup \cdots \cup R_\theta^{(i_{\ell})} \cup \overline{R},$$
with $0\leq \theta<\delta$, $\overline{R}\subset R_{\theta+1}^{(i_{\ell})}$.
Denote by $X$ and $Y$ the quantities $\sum_{i=1}^{\delta}r_i$ and $\sum_{i=1}^{\theta+1}r_i$ respectively. 
From (\ref{rankN}) we have

\begin{eqnarray*}
k-1&=&rank(\langle g_j\ : \ j \in \overline{N}\rangle)\leq \sum_{r=1}^{\ell-1} rank(\langle g_j\ : \ j \in N_{\delta}^{(i_r)}\rangle)+\sum_{r=1}^{\theta+1} rank(\langle g_j\ : \ j \in R_{r}^{(i_\ell)}\rangle)\\
&\leq& (\ell-1)\left (1+X\right)+Y,
\end{eqnarray*}
that is 
$$\ell-1 \geq \left \lceil \frac{k-1-Y}{1+X}\right\rceil.$$

Since $g_{i_r}\in \bigcap_{a=1}^{\delta}  \langle g_j\ : \ j \in R_a^{(i_r)}\rangle $, the difference $\#N_{a}^{i_r}-rank ( \langle g_j\ : \ j \in N_{a}^{i_r}\rangle)$ is at least $a$.

Therefore,
\begin{eqnarray*}
\#\overline{N}&\geq& rank(\langle g_j\ : \ j \in \overline{N}\rangle)+(\ell-1) \delta+\theta \\
&\geq &k-1+ \left \lceil \frac{k-1-Y}{1+X}\right\rceil \delta+\theta\\
&\geq &k-1+ \frac{k-1-Y}{1+X}\delta+\theta\\
&= &k-2+ \frac{(k-1)\delta+1}{1+X} + \frac{(1+X)(1+\theta)-1-Y\delta}{1+X}.
\end{eqnarray*}
Consider now $(1+X)(1+\theta)-1-Y\delta$. 
\begin{itemize}
\item Suppose first that $\theta=\delta-1$. Then $(1+X)(1+\theta)-1-Y\delta=(1+X)\delta-1-Y\delta \geq \delta-1$, since $X \geq Y$.
\item If $\delta \geq \theta+2$,
\begin{eqnarray*}
(1+X)(1+\theta)-1-Y\delta &=&X(\theta+1)- Y \delta+\theta\\
&=& (\theta+1) \sum_{i=1}^{\theta+1}r_i + (\theta+1)\sum_{i=\theta+2}^{\delta} r_i-\delta \sum_{i=1}^{\theta+1} r_i +\theta\\
&=& -(\delta-\theta-1)\sum_{i=1}^{\theta+1} r_i+ (\theta+1)\sum_{i=\theta+2}^{\delta} r_i+\theta\\
&\geq& -(\delta-\theta-1)(\theta+1)r_{\theta+1}+(\theta+1)(\delta-\theta-1)r_{\theta+2}+\theta\\
&\geq& \theta.
\end{eqnarray*}
\end{itemize}
So, $(1+X)(1+\theta)-1-Y\delta\geq 0$.
Since $\#\overline{N}$ is an integer,
\begin{eqnarray*}
\#\overline{N}&\geq& k-2+ \left \lceil \frac{(k-1)\delta+1}{1+X} \right \rceil.
\end{eqnarray*}
By Lemma \ref{Lemma:MinimumDistance},
\begin{eqnarray*}
n &=& d + max\{\#N \  : \ N\subset\{1,\ldots,n\}, \ rank(\langle g_j\ : \ j \in N\rangle) < k)\geq d+\#\overline{N}\\
&\geq &d+k-2+\left \lceil \frac{(k-1)\delta+1}{1+X}\right\rceil,
\end{eqnarray*}
and  the claim follows.
\endproof

\begin{remark}
Note that for $\delta=1$, Bound \eqref{Eq:MainBound} is weaker than Bound \eqref{rec1}.
\end{remark}

In \cite{JKZ2019} the authors construct LRC codes over the rational function field $\fq(x)$ (genus zero) dealing with subgroup of its automorphism group, and obtain codes with length $n \simeq q$ roughly and relative defect  $\simeq 1/q$ according to Formula \eqref{Eq:RSD}.

\section{General construction from subgroups with trivial intersection}

Many constructions of LRC codes from function fields arose in recent years, over the rational function field \cite{JKZ2019, JKZ}, elliptic function fields \cite{LLX2019} and over algebraic curves (see \cite{MT2018, WZ2014}), fiber product of curves \cite{HMM2018} and curves with separated variables \cite{MTT2019}. In a variant of the different ways to construct LRC codes, we focus on codes from algebraic function fields of genus $g \geq 1$ using certain subgroups of the automorphism group of the underlying curve. In this section we deal with the case of subgroups with trivial intersection and in the next Section with non-trivial intersection subgroups.
\\

%%%%%%%%%%%%%%%%%%%%%%%%%%%%%%%%%%%%%%%%%%%%%%%%%%%
%%%%%%%%	MAIN RESULT-TRIVIAL INTERSECTION	 	%%%%%%%%%%%%%%%%
%%%%%%%%%%%%%%%%%%%%%%%%%%%%%%%%%%%%%%%%%%%%%%%%%%%
\begin{theorem}\label{Th:GeneralConstruction}
Let  $\cF | \fq$ be a function field of genus $g$.  Consider $s$ subgroups $H_i$ of the automorphism group of $\cF| \fq$, each of sizes $r_i+1$, such that the group $\cG \simeq \bigotimes_{i=1}^s\cH_i$ is isomorphic to the internal direct product of $\cH_1,\ldots,\cH_s$.  Let $\mathcal{P}$ be a set of places in $\cF$ lying over $m$ rational places in the fixed field $\cF^{\cG}$ that are completely split in the extension $\cF|\cF^{\cG}$. Define $n=m \prod _{i=1}^s (r_i+1)$, that is, $n$ to be the total number of places of $\cF$ lying over the $m$ selected rational places in $\cF^{\cG}$.

Suppose that there exists a place $P_{\infty}$ of $\cF$ which is completely ramified in $\cF |\cF^{\cG}$ and let $Q_{\infty}^{(i)}$ be the unique place in $\cF^{\cH_i} $ lying under $P_{\infty}$. For $i=1,\ldots, s$ suppose further there exist functions  $z_i, w_i$, such that 
\begin{enumerate}[i)]
\item $z_i\in \cF^{\cH_i}$, $supp((z_i)_{\infty})=\{Q_{\infty}^{(i)}\}$;
\item $w_i\in \cF \setminus \cF^{\cH_i}$, $supp((w_i)_{\infty})=\{P_{\infty}\}$;
\item $w_i : P^{\cH_i} \to \mathbb{F}_q$ is injective.
\end{enumerate}
Let $t_i \geq 1$ be such that  
$$V_i := \left\{ \sum_{\ell=0}^{r_i-1} \left( \sum_{j=0}^{t_i} a_{\ell j}^{(i)}z_i^j \right)w_i^\ell \in \cF \mid  a_{\ell j}^{(i)} \in \fq \right\}\subset \mathcal{L}((n-d)P_{\infty})$$
for some $1 \leq d \leq n$ and let $V=\bigcap_{i=1}^s V_i \subset \mathcal{L}((n-d)P_{\infty})$.
If $\dim_{\mathbb{F}_q}(V)>0$, then there exists an 
$$[n,\dim_{\mathbb{F}_q}(V),\geq d;r_1,\ldots,r_s]\text{-recoverable code.}$$
\end{theorem}

\proof
Let $V=\bigcap_{i=1}^s V_i \subset \mathcal{L}((n-d)P_{\infty})$ and consider the following linear map
$$\begin{array}{llll}
e_{\mathcal{P}} \ : \ & V &\to& \mathbb{F}_q^n\\
&f &\mapsto &e_{\mathcal{P}}(f)=(f(P_1),\ldots, f(P_n)),
\end{array}$$
where $\mathcal{P}=\{P_1,\ldots, P_n\}$. 

The linear code $e_{\mathcal{P}}(V)$ is contained in the algebraic geometry code $C_\cL(\sum_i^n P_i, (n-d)P_{\infty})$, so by the bound on Equation \eqref{goppabound} it has minimum distance at least $d$. By definition it has length $n$ and dimension $\dim_{\mathbb{F}_q}(V)$. 

Now, we deal with recoverability. First note that a rational place $Q\in\cF^{\cG}$ which is completely split in $\cF|\cF^{\cG}$ satisfies $P^{\cH_i}\cap P^{\cH_j}=\{P\}$ for each $i\neq j \in \{1,\ldots, s\}$, where $P$ is an arbitrary place in  $\cP$ lying over $Q$. In fact, since $Q$ is completely split, by the orbit stabilizer theorem, this means that the stabilizer in $G$ of  $P$ is trivial. If $P^{H_i}$ and $P^{H_j}$ contain another common place $\tilde P \ne P$ then there exist  $h_i \in H_i$ and $h_j \in H_j$ such that $\tilde P=h_i(P)=h_j(P)$ so that $h_i h_j^{-1}(P)=P$. Since $H_i \cap H_j$ is trivial, $h_i h_j^{-1}$ is a non-trivial element of $G$, a contradiction.
Consider a fixed place $P \in \mathcal{P}$ and an $f \in V$. We are going to show that $f(P)$ can be repaired by any of the $s$ recoverable sets (depending on $P$)
$$R_i = \{f(P^{\prime}) \ : \ P^{\prime} \in P^{\cH_i} \setminus\{P\}\},\text{ for } $$  
of cardinality $r_i$ for $i=1, \dots, s$.
Since $f \in V_i$ we can write $f=\sum_{\ell=0}^{r_i-1} \left( \sum_{j=0}^{t_i} a_{\ell j}^{(i)}z_i^j \right)w_i^\ell, a_{\ell j}^{(i)} \in \fq$. Since $z_i\in \cF^{\cH_i}$ we have that $z_i(Q_1)=z_j(Q_2)=\overline{z}$ for any $Q_i, Q_j \in P^{\cH_i}$, then for each $Q \in P^{\cH_i}$ we have
$$f(Q)=\sum_{\ell=0}^{r_i-1} \left( \sum_{j=0}^{t_i} a_{\ell j}^{(i)}z_i^j (Q)\right)w_i(Q)^\ell =\sum_{\ell=0}^{r_i-1} \left( \sum_{j=0}^{t_i} a_{\ell j}^{(i)}\overline{z}^j\right)w_i(Q)^\ell=\sum_{\ell=0}^{r_i-1} \gamma_{\ell,P}^{(i)} w_i(Q)^\ell,$$
where $\gamma_{\ell,P}^{(i)}\in \mathbb{F}_q$ depends only on $\ell$ (and $P$ clearly). So, 
$$h(X)=\sum_{\ell=0}^{r_i-1} \gamma_{\ell,P}^{(i)} X^\ell\in \mathbb{F}_q[X]$$
is a polynomial of degree at most $r_i-1$. Since the values  $w_i(Q)$ with  $Q\in  P^{\cH_i}$ are distinct we know exactly $r_i$ different values of $h(X)$. By Lagrange interpolation we can recover $h(X)$ and then compute the missing value $h(w_i(P))$. By definition of $\mathcal{P}$, the recovering sets $  P^\mathcal{H}_i \setminus \{P\}$ are pairwise disjoint.
\endproof

In what follows we are going to construct examples of LRC codes over algebraic function fields of genus $g \geq 1$ based on Theorem \ref{Th:GeneralConstruction}. 
All the constructed codes comes from algebraic curves with many rational points.

%%%%%%%%%%%%%%%%%%%%%%%%%%%%%%%%%%%%%%%%%%%%%%%%%%%
%%%%%%%%%%%%		HERMITIAN 1	 	%%%%%%%%%%%%%%%%%%%%%%%%%
%%%%%%%%%%%%%%%%%%%%%%%%%%%%%%%%%%%%%%%%%%%%%%%%%%%
%\begin{example}\label{Ex:Hermitian1}

\subsection{Codes from the Hermitian Function Field I}\label{Ex:Hermitian1}

Consider the Hermitian function field $\cF=\mathbb{F}_{q^2}(x,y)$ with $y^{q+1}=x^q+x$. Choose a divisor $u\geq 2$ of $q+1$. Consider  the following two subgroups of $\fq$-automorphisms of $\cF$
\begin{eqnarray*}
    	\cH_1&:=&\{(x,y) \mapsto (x+c,y) \mid  c \in \mathbb{F}_{q^2}, c^q+c=0 \},\\
    	\cH_2&:=&\{(x,y) \mapsto (x,a^i y)\mid  i=0,\ldots,u-1 \text{ and } a \in \fqs^* \text{ of order } u \geq 2\}.
\end{eqnarray*}
It is easy to see that $\# \cH_1=q$, $\#\cH_2=u$, $\mathcal{F}^{\cH_1}=\mathbb{F}_{q^2}(y)$ and $\mathcal{F}^{\cH_2}=\mathbb{F}_{q^2}(x,y^u)$.  The groups have trivial intersection and commute, then $\cG=\cH_1\cH_2=\cH_1 \times \cH_2$  and $\cF^\cG=\fqs(z:=y^u) \subseteq \cF^{\cH_i},$ for $ i=1, 2$.  Let $P_{\infty}$ be the unique pole of $x$ and $y \in \cF$, then $(y)_\infty = qP_\infty $ and  $(x)_\infty=(q+1)P_\infty $ in $\cF$. Denote by $Q_\infty^i=P_\infty \cap \cF^{\cH_i},$ for $ i=1, 2$.  In the extension $\cF \mid\fqs(z)$, we have that $(q^2-1)/u$ places in $\cF^\cG$ are totally split, so $\#\mathcal{P}=q(q^2-1)=n$. For $t_i \geq 2, i=1, 2$ consider the $\fq$-vector spaces
\begin{eqnarray*}
	V_1 &:=&\left\{\sum_{\ell=0}^{q-2} \left(\sum_{j=0}^{t_1}a_{\ell,j}y^j  \right)x^\ell   \ : \ a_{\ell,j} \in \mathbb{F}_{q^2}\right \} \text{ and } \\
	V_2 &:=&\left\{\sum_{\ell=0}^{u-2} \left(\sum_{j=0}^{t_2}a_{\ell,j}z^{j}  \right)y^\ell   \ : \ a_{\ell,j} \in \mathbb{F}_{q^2}\right \}=\left\{\sum_{\ell=0}^{u-2} \sum_{j=0}^{t_2} a_{\ell,j}y^{(u-1)j+\ell}   \ : \ a_{\ell,j} \in \mathbb{F}_{q^2}\right\}.\\
\end{eqnarray*}
Note that $\dim_{\fqs} V_1=(q-1)(t_1+1) $ and $\dim_{\fqs} V_2=(u-1)(t_2+1)$.
A function in $V_1$ has pole only at the place $P_\infty$ of order at most $q^2-q-2+t_1q$, while a function in $V_2$ also has pole only at $P_\infty$ of order at most $q((u-1)t_2+u-2)$. 

Using the notations of the Theorem \ref{Th:GeneralConstruction}, we have $z_1=y$, $z_2=y^{u-1}$, $w_1=x$ and $w_2=y$. Now we choose the parameters:
\begin{enumerate}[i)]
\item $t_1=q-1 \, \text{ and } \, \frac{2q^2-qu-2}{q(u-1)}\leq t_2 < \frac{q^2}{u-1}-1$
so that $$V_1,V_2 \subset \mathcal{L}(q((u-1)t_2+u-2))P_{\infty}), \text{ and } d=n-q(ut_2+u-t_2-2)\geq 1.$$
Now, 
$$V=V_1\cap V_2=\left\{\sum_{j=0}^{q-1}a_{j}y^j  \ : \ a_{j} \in \mathbb{F}_{q^2}\right \} \textrm{ and } \dim _{\mathbb{F}_{q^2}}V=q.$$
\item $u(t_2+1)=q+1, t_1 \geq 1$, $(t_1,t_2)\neq (1,1)$, then 
$$V_1,V_2 \subset \mathcal{L}((q^2-q-2+t_1q)P_{\infty}), $$
 and $V=V_1\cap V_2=\left\{\sum_{j=0}^{q-1}a_{j}y^j  \ : \ a_{j} \in \mathbb{F}_{q^2}\right \}$. 
\end{enumerate}
The two choices of parameters in the construction can be summarized as the following.
\begin{proposition}\label{herm1}
For every $1<u$ a divisor of $q+1$ and $\frac{2q^2-qu-2}{q(u-1)}\leq t_2 < \frac{q^2}{u-1}-1$ there exists a 
$$[q(q^2-1), q, \geq q(q^2-1)-q(ut_2+u-t_2-2); q-1, u-1]\text{-recoverable code }\mathcal{C}_1.$$  
For every $1<u, u(t_2+1)=q+1, t_1 \geq 2$ and $(t_1,t_2)\neq (1,1)$ there exists a 
$$[q(q^2-1), q, \geq q(q^2-1)-(q^2-q-2+t_1q); q-1, u-1]\text{-recoverable code }\mathcal{C}_2.$$
\end{proposition}

\begin{remark}The relative defect of the codes $\cC_1$ and $\cC_2$ (see Equation \eqref{Eq:RSD}) are 
\begin{align*}
\Delta(\cC_1)&\leq \frac{1}{q(q^2-1)}\left(q(ut_2+u-t_2-2)-q+2-\left\lceil\frac{2q-1}{q+u-1}\right\rceil\right)=\frac{q(ut_2+u-t_2-3)+1}{q(q^2-1)},\\
&\\
\Delta(\cC_2)&\leq \frac{1}{q(q^2-1)}\left((q^2-q-2+t_1q)-q+2-\left\lceil\frac{2q-1}{q+u-1}\right\rceil\right)=\frac{q^2-2q+t_1q-1}{q(q^2-1)}
\simeq \frac{q+t_1}{q^2-1}.
\end{align*}

For the code $\cC_1$, if $t_2$ is close to the lower bound $\frac{2q^2-qu-2}{q(u-1)}$, $\delta(\mathcal{C}_1)$ is less or equal to $2/q$. 
\end{remark}

%%%%%%%%%%%%%%%%%%%%%%%%%%%%%%%%%%%%%%%%%%%%%%%%%%%
%%%%%%%%%%%%		HERMITIAN 2	 	%%%%%%%%%%%%%%%%%%%%%%%%%
%%%%%%%%%%%%%%%%%%%%%%%%%%%%%%%%%%%%%%%%%%%%%%%%%%%
\subsection{Codes from the Hermitian Function Field II}\label{Ex:Hermitian2}

Consider the Hermitian function field $\cF=\mathbb{F}_{q^2}(x, y)$ with $y^{q+1}=x^q+x$. Use the same notations as in Section \ref{Ex:Hermitian1}. Consider now
\begin{eqnarray*}
	V_1 &:=&\left\{\sum_{\ell=0}^{q-1} \left(\sum_{j=0}^{t_1 \leq u}a_{\ell,j}y^j  \right)x^\ell   \ : \ a_{\ell,j} \in \mathbb{F}_{q^2}\right \},\\
	V_2 &:=&\left\{\sum_{\ell=0}^{u-1} \left(\sum_{j=0}^{t_2 \leq q}a_{\ell,j}x^j  \right)y^\ell   \ : \ a_{\ell,j} \in \mathbb{F}_{q^2}\right \},\\
\end{eqnarray*}
of dimensions $\dim _{\mathbb{F}_{q^2}} V_1= (t_1+1)q$ and $\dim _{\mathbb{F}_{q^2}} V_2= (t_2+1)u$. 

A function in $V_1$ has pole only at the place $P_\infty$ of order at most $q^2-q-2+t_1q$, while a function in $V_2$ also has pole only at $P_\infty$ of order at most $q((u-1)t_2+u-2)$.

Using the notations of the Theorem \ref{Th:GeneralConstruction} we have $z_1=y$, $z_2=x$, $w_1=x$, $w_2=y$.  Also, 
$$V=V_1\cap V_2=\left\{\sum_{\ell=0}^{t_1} \left(\sum_{j=0}^{t_2}a_{\ell,j}x^j  \right)y^\ell   \ : \ a_{\ell,j} \in \mathbb{F}_{q^2}\right \}$$
and $\dim _{\mathbb{F}_{q^2}} V=(t_1+1)(t_2+1)$. For  
$$M=\max\{t_1 q+q^2-1,(u-1)q+t_2(q+1)\}\leq 2q^2+q-1$$
we have $V \subset \mathcal{L}(MP_{\infty}).$
In this case $n-d=M$ and  $ d \geq \max\{q^2+t_1(q-t_2),(u-1)q+t_2(q-t_1)\}$.
We obtain the following proposition.
\begin{proposition}\label{herm2}
For every divisor $1<u$ of $q+1$, $1\leq t_1 \leq u$, $1\leq t_2 \leq q$ and $$M=\max\{t_1 q+q^2-1,(u-1)q+t_2(q+1)\}$$ there exists a 
$$[q(q^2-1), (t_1+1)(t_2+1), \geq q(q^2-1)-M; q-1, u-1]\text{-recoverable code }\mathcal{C}.$$
\end{proposition}
\begin{remark}
The relative defect of this code satisfies
$$\Delta(\mathcal{C})\leq \frac{M-(t_1+1)(t_2+1)+2-\left\lceil \frac{2(t_1t_2+t_1+t_2)+1}{q+u-1}\right\rceil}{q(q^2-1)}.$$
Choosing $u=q+1, t_1=u$ and $t_2=q$ we obtain $\Delta(\mathcal{C})\leq \frac{q^2-q-1}{q(q^2-1)}$ and so the relative defect goes to zero as $q$ goes to infinity.
\end{remark}

%%%%%%%%%%%%%%%%%%%%%%%%%%%%%%%%%%%%%%%%%%%%%%%%%%%
%%%%%%%%%%%%		GK 1			 	%%%%%%%%%%%%%%%%%%%%%%%%%
%%%%%%%%%%%%%%%%%%%%%%%%%%%%%%%%%%%%%%%%%%%%%%%%%%%
\subsection{Codes from the Giulleti-Korchmáros curve}\label{Ex:GK1}
Let $K= \mathbb{F}_{q^{6}}$ and $\cF=K(x, y, z)$ be the function field of the curve $\mathcal{GK}$ whose affine model is given by the complete intersection
$$
\left\{ \begin{array}{l}
Z^{q^2-q+1}=Y^{q^2}-Y,\\
Y^{q+1}=X^q+X.
\end{array}
\right. $$
This is a maximal curve over $\mathbb{F}_{q^{6}}$ and has $q^8-q^6+q^5+1$ rational places and only one place $P_\infty$ at infinity. Moreover, $P_\infty$ is the common pole of $x, y $ and $z$ with pole divisors 
$$(x)_{\infty}=(q^3+1)P_{\infty}, \qquad  (y)_{\infty}=(q^3-q^2+q)P_{\infty}, \qquad (z)_{\infty}=qP_{\infty}, $$see \cite{GK2009}. 
In this case we consider three subgroups of the following types, let $A$ be a subgroup of $\{a \in \mathbb{F}_{q^6} \ : \ a^q+a=0\}$,  $\eta, \omega \in \mathbb{F}_{q^6}$, with $ord(\eta) \mid q^3+1$,  $ord(\omega) \mid q^2-q+1$, and $\gcd(ord(\eta), ord(\omega))=1$. Consider the following subgroups of $Aut(\cF | K)$: 
\begin{eqnarray*}
	\cH_1&:=&\{\sigma_a: (x, y, z) \mapsto (x+a, y, z) \ : \ a \in A\};\\
	\cH_2&:=&\{\sigma_{i}: (x, y, z) \mapsto (x,\eta^{i(q^2-q+1)} y,\eta^{i}z)\ : \  i=0,\ldots, ord(\eta)-1\};\\
	\cH_3&:=&\{\sigma_{i}: (x, y, z) \mapsto (x, y,\omega^i z)\ : \ i=0,\ldots, ord(\omega)-1\},
\end{eqnarray*}
Then we have $\cG=\cH_1\cH_2\cH_3 \cong \cH_1\times \cH_2\times \cH_3$ is a subgroup of $Aut(\cF | K)$. The fixed fields satisfies $K(y, z) = \cF^{\cH_1}, K(x) \subseteq \cF^{\cH_2}$ and $ K(x, y) \subseteq \cF^{\cH_3}$. 
Now we notice that for all $P=(x, y, z)$ a rational place with $z \ne 0$ we have that $P^{\cH_i} \cap P^{\cH_j}=\{P\}, i, j \in \{1, 2, 3\}$. In fact, if $Q \in P^{\cH_2} \cap P^{\cH_3}$, since the unique possibility for  $\eta^i=1$ and $\omega^j=1$ yields $i=j=0$ from $\gcd( ord(\eta), ord(\omega))=1$, we obtain $Q=P$. The other cases are trivial. Therefore, all the affine $\mathbb{F}_{q^6}$-rational points of  $\mathcal{GK}$ belong to long orbits with respect to $\cG$ are those with $z\neq 0$. So, the cardinality of the set  $\cP$ as in Theorem \ref{Th:GeneralConstruction} can be taken as $n=q^8-q^6+q^5-q^3=q^3(q^2-1)(q^3+1)$.  
Using the notations in Theorem \ref{Th:GeneralConstruction} we consider 
 \begin{enumerate}[i)]
 \item $z_1 :=z \in \cF^{\cH_1}, w_1:= x \in \cF \setminus \cF^{\cH_1}$
 \item $z_2 := x \in \cF^{\cH_2}, w_2:= z \in \cF \setminus \cF^{\cH_2}$
 \item $z_3 :=x \in \cF^{\cH_3}, w_1:= z \in \cF \setminus \cF^{\cH_3}$.
\end{enumerate}
For any integers $t_i \geq 0, i=1, 2, 3$, consider the following $\mathbb{F}_{q^{6}}$-vector spaces
\begin{eqnarray*}
	V_1 &:=&\left\{\sum_{\ell=0}^{\#A-2} \left(\sum_{j=0}^{t_1}a_{\ell, j}z^j  \right)x^\ell   \ : \ a_{\ell, j} \in \mathbb{F}_{q^6}\right \},\\
	V_2 &:=&\left\{\sum_{\ell=0}^{ord(\eta)-2} \left(\sum_{j=0}^{t_2}a_{\ell, j}x^j  \right)z^\ell   \ : \ a_{\ell, j} \in \mathbb{F}_{q^6}\right \}, \text{ and }\\
	V_3 &:=&\left\{\sum_{\ell=0}^{ord(\omega)-2} \left(\sum_{j=0}^{t_3}a_{\ell, j}x^j  \right)z^\ell   \ : \ a_{\ell, j} \in \mathbb{F}_{q^6}\right \}.\\
\end{eqnarray*}
For $N_1=\min\{\#A-2, t_2, t_3\}$, $M_1=\min\{t_1, ord(\eta)-2, ord(\omega)-2\}$ and  $$S = \max\{(\#A-2)(q^3+1)+ t_1q, t_2(q^3+1)+(ord(\eta)-2)q, t_3(q^3+1)+(ord(\omega)-2)q\}$$ we obtain
$$V:=V_1\cap V_2 \cap V_3 =\left\{\sum_{n_1=0}^{N_1} \sum_{m_1=0}^{M_1}a_{n_1,m_1}z^{m_1}  x^{n_1}  \ : \ a_{n_1,m_1} \in \mathbb{F}_{q^6}   \right\}  \subseteq \cL( SP_\infty)$$
has dimension $\dim_{\mathbb{F}_{q^6}}V=(M_1+1)(N_1+1)$. 
Choosing $t_i, i=1, 2, 3$ such that $q^8-q^6+q^5-q^3-S \geq 1$ we obtain the following proposition.
\begin{proposition}\label{GK}
Let $q=p^\ell, p $ prime. Then for any 
\begin{enumerate}[i)]
\item $a=p^h, 1\leq h \leq \ell$, 
\item $ord(\eta) \mid q^3+1$, $ord(\omega) \mid q^2-q+1$, with $\gcd(ord(\eta), ord(\omega))=1$, 
\item $0< t_1 < q^2(q^2-1)(q^3+1)-q^3+2q^2-1,  0< t_2 < q^3(q^2-1)-q,  0< t_3 <q^3(q^2-1)-1$,
\item $N_1=\min\{a-2, t_2, t_3\}$, $M_1=\min\{t_1, ord(\eta)-2, ord(\omega)-2\},$ and
\item $S = \max\{(a-2)(q^3+1)+ t_1q, t_2(q^3+1)+(ord(\eta)-2)q, t_3(q^3+1)+(ord(\omega)-2)q\}$
\end{enumerate}
there exists a 
$$[n:=q^8-q^6+q^5-q^3, (M_1+1)(N_1+1), \geq n-S; a-1, ord(\eta)-1, ord(\omega)-1]$$
recoverable code $\mathcal{C}_3$ over $\mathbb{F}_{q^6}$.
\end{proposition}

Its relative defect is 
$$\Delta(\mathcal{C}_3)\leq \frac{S-(M_1+1)(N_1+1)+2-\left\lceil \frac{3(M_1N_1+M_1+N_1)+1}{a+ord(\eta)+ord(\omega)-2}\right\rceil}{q^8-q^6+q^5-q^3}.$$
Suppose that $(q^3+1)$ possesses two coprime divisors both close to $q\sqrt{q}$. Consider $q\sqrt{q}<t_1<q^2$, $t_2, t_3\simeq a$, so that $N_1=a-2$, $M_1\simeq q\sqrt{q}$, $S\simeq N_1q^3$ and therefore $\Delta(\mathcal{C}_3)$ is at most a value close to $N_1(q^3-q\sqrt{q})/q^8.$
%
%\textcolor{blue}{As a particular case of the proposition above we consider the following.}
%
%\textcolor{blue}{
%\begin{proposition}\label{GKPrime}
%Let $q=p^\ell, p $ prime. Then for any 
%\begin{enumerate}[i)]
%\item $ord(\eta) \mid q^3+1$, $ord(\omega) \mid q^2-q+1$, with $\gcd(ord(\eta), ord(\omega))=1$, 
%\item $0< t_1 < q^2(q^2-1)(q^3+1)-q^3+2q^2-1,  0< t_2 < q^3(q^2-1)-q,  0< t_3 <q^3(q^2-1)-1$,
%\item $N_2=\min\{ t_2, t_3\}$, $M_2=\min\{ord(\eta)-2, ord(\omega)-2\},$ and
%\item $S = \max\{t_2(q^3+1)+(ord(\eta)-2)q, t_3(q^3+1)+(ord(\omega)-2)q\}$
%\end{enumerate}
%there exists a 
%$$[n:=q^8-q^6+q^5-q^3, (M_2+1)(N_2+1), \geq n-S; ord(\eta)-1, ord(\omega)-1)]$$
%recoverable code $\mathcal{C}_3^{\prime}$ over $\mathbb{F}_{q^6}$.
%\end{proposition}
%}
%\textcolor{blue}{
%\proof
%In this case we consider $V=V_2\cap V_3$ (see the notation above). That is,
%$$V=\left\{\sum_{n_2=0}^{N_2} \sum_{m_2=0}^{M_2}a_{n_2,m_2}z^{m_2}  x^{n_2}  \ : \ a_{n_2,m_2} \in \mathbb{F}_{q^6}   \right\}  \subseteq \cL( SP_\infty),$$
%where $N_2=\min\{t_2, t_3\}$, $M_2=\min\{ord(\eta)-2, ord(\omega)-2\}$.
%So,  $\dim_{\mathbb{F}_{q^6}}V=(M_2+1)(N_2+1)$ and the minimum distance  of  $\mathcal{C}_3^{\prime}$ is at least $ n-S=q^8-q^6+q^5-q^3-S$, where   $$S = \max\{ t_2(q^3+1)+(ord(\eta)-2)q, t_3(q^3+1)+(ord(\omega)-2)q\}.$$
%Choosing $t_i, i=1, 2, 3$, such that $q^8-q^6+q^5-q^3-S \geq 1$, the claim follows.
%\endproof}

\begin{remark}\label{Remark:Comparison}
In \cite{HMM2018} the authors consider fiber products of algebraic curves in order to construct LRC codes. In particular \cite[Corollary 1]{HMM2018} can be seen as a construction from trivially intersecting subgroups of the automorphism group of an algebraic curve. In what follows we choose codes with same length and compare the obtained relative defects. 

In \cite[Theorem 5.1]{HMM2018} the authors construct codes on the generalized GK curve. For the special case of the curve GK with $\ell=q$, they obtain a 
$$
[n:=q^8-q^6+q^5-q^3, (q^2-q)(q^2-1), n-q^5+2q^3-q^2-q+2; q-1, q^2-q]
$$
recoverable code over $\mathbb{F}_{q^6}$. 
%with relative singleton defect less or equal than
%$ \frac{q^5-q^4-q^3-3q-1}{q^8-q^6+q^5-q^3}.$
The relative defect of such a code is $\frac{q^5-q^4-q^3+2q+1}{n}$.

There are many choices for $a,\eta,\omega$ such that the code constructed from Proposition \ref{GK} has relative defect smaller than $\frac{q^5-q^4-q^3+2q+1}{n}$. Suppose $3\nmid (q+1)$. Consider for instance $a=q$, $ord(\eta)=q+1$, $ord(\omega)=q^2-q+1$, $t_1=t_2=t_3=q-1$. In this case, $N_1=q-2$, $M_1=q-1$ and $S=q^4-q^2-1$. Note that $S$ is the maximum $S=\max\{(a-2)(q^3+1)+ t_1q, t_2(q^3+1)+(ord(\eta)-2)q, t_3(q^3+1)+(ord(\omega)-2)q\}=t_3(q^3+1)+(ord(\omega)-2)q=(q-1)(q^3+1)+(q^2-q-1)q=q^4-q^2-1$, and hence there exists an
$$[n:=q^8-q^6+q^5-q^3, q(q-1), \geq q^8-q^6+q^5-q^4-q^3+q^2+1; q-1, q, q^2-q]$$
recoverable code $\mathcal{C}_3$ over $\mathbb{F}_{q^6}$ whose relative defect is at most 
$$\frac{S-(M_1+1)(N_1+1)+2-\left\lceil \frac{3(M_1N_1+M_1+N_1)+1}{a+ord(\eta)+ord(\omega)-2}\right\rceil}{n}\leq \frac{q^4-2q^2+q}{n}.$$
Similar results can be obtained when $3\mid (q+1)$ by considering for instance $ord(\eta)=q+1$ and $ord(\omega)=(q^2-q+1)/3$. Also in this case, the relative defect is roughly $q^4/n$.
\end{remark}

%%%%%%%%%%%%%%%%%%%%%%%%%%%%%%%%%%%%%%%%%%%%%%%%%%%
%%%%%%%%%%%%		KUMMER OF LINEARIZED POLYNOMIALS 	%%%%%%%%%%%%
%%%%%%%%%%%%%%%%%%%%%%%%%%%%%%%%%%%%%%%%%%%%%%%%%%%
%\subsection{Codes from Kummer extensions}\label{Ex:KummerLinerarized}
%
%\lu{since we do not know the number of rational points for a general curve of this type, I suggest to delete this example}
%
%Let $q=p^\ell$ and consider the function field $\cS$ of the algebraic curve defined by the affine equation $S: y^m=L(x)$, where  \lu{$m\mid q+1$} and  $L(x)$ is a separable $p$-linearized polynomial over $\fq$ of degree $p^h < q$ with all its roots in $\fq$.  

%This curve has genus $g=\frac{(m-1)(p^h-1)}{2}$ and more results on these type of curves can be found in \cite[Chapter 12]{HKT}. 
%We define two subgroups of the automorphism group of the curve by
%$$H_1=\{(x,y) \mapsto (x+a,y) \mid L(a)=0 \ \text{ and }  a \in \fq \},$$
%and
%$$H_2=\{(x,y) \mapsto (x,\lambda y) \mid \lambda \in \fq \text{ and } \lambda^m=1 \} ,$$
%so that $\# H_1=p^h$, $\# H_2=m$, $H_1H_2=H_2H_1$ and they have trivial intersection. It is easy to see that the fixed fields are $\cS^{H_1}=\fq(y)$ and $\cS^{H_2}=\fq(x)$. Let $P_\infty $ be the only pole of $x$ and $y$ in the function field $\mathcal{S}$, this place is fixed by $H_1$ and $H_2$ and its Weierstrass semigroup $H(P_\infty)=\langle p^h, m \rangle$  considering the divisors of $x$ and $y$. 

%%%%%%%%%%%%%%%%%%%%%%%%%%%%%%%%%%%%%%%%%%%%%%%%%%%
%%%%%%%%%%%%		GENERALIZED HERMITIAN 	%%%%%%%%%%%%%%%%%%%
%%%%%%%%%%%%%%%%%%%%%%%%%%%%%%%%%%%%%%%%%%%%%%%%%%%
\subsection{Codes from the generalized Hermitian curve}\label{Ex:GenHermitian}
Let $q$ be odd. The curve $S: y^{q^\ell+1}=x^q+x$ over $\mathbb{F}_{q^{2\ell}} $ with $\ell \geq 1$ odd has $q^{2\ell +1}+1$ rational points and genus $g=q^\ell(q-1)/2$, see \cite{KKO2001, HKK2006}. This curve has only one rational point at infinity denoted by $P_\infty$, it is also the only pole of the functions $x$ and $y$ with pole order 
$q^\ell+1$ and $q$, respectively.
Let $\cS$ be the function field of the curve $S$, and consider two subgroups of the automorphism group of the curve given by
\begin{eqnarray*}
H_1&:=&\{(x, y) \mapsto (x+a, y) \mid a^q+a=0 \ \text{ and }  a \in \mathbb{F}_{q^{2\ell}} \}, \text{ and }\\
H_2&:=&\{(x, y) \mapsto (x, \lambda y) \mid \lambda \in \mathbb{F}_{q^{2\ell}} \text{ and } \lambda^{q^\ell+1}=1 \} ,
\end{eqnarray*}
so that  $\# \cH_1=q, \# \cH_2=q^\ell+1$ and the fixed fields are $\cS^{H_1}=\mathbb{F}_{q^{2\ell}}(y)$, $\cS^{H_2}=\mathbb{F}_{q^{2\ell}}(x)$, and $S^\cG=\mathbb{F}_{q^{2\ell}}(x^q+x)=\mathbb{F}_{q^{2\ell}}(y^{q^\ell+1})$ by the direct product $\cG \simeq \cH_1 \times \cH_2$. The number of rational places in  $S^\cG$ that completely splits in the field extension $S| S^\cG$ is $q^\ell-1$, then we have $n=\# \cP=q^{2\ell +1}-q$.
Let  $ 0 \leq t_i, i=1, 2$ and consider the vector spaces
\begin{eqnarray*}
	V_1 &:=&\left\{\sum_{\ell=0}^{q-2} \left(\sum_{j=0}^{t_1 }a_{\ell, j}y^j  \right)x^\ell   \ : \ a_{\ell, j} \in \mathbb{F}_{q^{2\ell}}\right \},\\
	V_2 &:=&\left\{\sum_{\ell=0}^{q^\ell-1} \left(\sum_{j=0}^{t_2 }a_{\ell, j}x^j  \right)y^\ell   \ : \ a_{\ell, j} \in \mathbb{F}_{q^{2\ell}}\right \},\\
\end{eqnarray*}
 Choose $t_i$ such that for $M_1=\max\{t_1, q^\ell-1 \}, M_2=\max\{t_2, q-2 \}$  we have $S=M_1q+M_2(q^\ell-1) \leq q^{2\ell +1}-q$. Then  $V_i \subset \cL(SP_\infty),$ for $ i=1, 2$, and for  $m_1=\min\{t_1, q^\ell-1 \}, m_2=\min\{t_2, q-2 \}$ we obtain
$$V=V_1\cap V_2=\left\{\sum_{\ell=0}^{m_1} \left(\sum_{j=0}^{m_2 }a_{\ell, j}x^j  \right)y^\ell   \ : \ a_{\ell, j} \in \mathbb{F}_{q^{2\ell}}\right \}  \textrm{ and } \dim _{\mathbb{F}_{q^2}}V=(m_1+1)(m_2+1).$$
The minimum distance of the code satisfies $d \geq q^{2\ell +1}-q- S$. We have obtained the following proposition.

\begin{proposition}\label{genherm}
Let $q$ be odd. Consider $0 \leq t_i, i=1, 2$ satisfying that 
$$S=M_1q+M_2(q^\ell-1) \leq q^{2\ell +1}-q$$ 
for $M_1=\max\{t_1, q^\ell-1 \}$ and $ M_2=\max\{t_2, q-2 \}$. Fixing $m_1=\min\{t_1, q^\ell-1 \}$ and $ m_2=\min\{t_2, q-2 \}$ there exists a 
$$[q^{2\ell +1}-q, (m_1+1)(n_1+1), \geq q^{2\ell +1}-q-S; q-1 , q^\ell]$$
recoverable code $\mathcal{C}$ over $\mathbb{F}_{q^{2\ell}}$.
\end{proposition}

The codes arising from Proposition \ref{genherm} have relative designed distance
$$\Delta(\mathcal{C})\leq \frac{S-(m_1+1)(n_1+1)+2-\left\lceil \frac{2(m_1n_1+m_1+n_1)+1}{q^\ell+q}\right\rceil}{q^{2\ell +1}-q}.$$

%%%%%%%%%%%%%%%%%%%%%%%%%%%%%%%%%%%%%%%%%%%%%%%%%%%
%%%%%%%%%%%%		ELLIPTIC CURVES 	%%%%%%%%%%%%%%%%%%%%%%%%%
%%%%%%%%%%%%%%%%%%%%%%%%%%%%%%%%%%%%%%%%%%%%%%%%%%%
%\subsection{Codes from Elliptic curves}\label{Ex:EllipticCurves}

%\lu{there is a full paper on optimal LRC codes over elliptic curves, I suggest we leave this out}
%\textbf{Elliptic curves with zero $j$-invariant.} Let $p \ne 2,3$. In \cite[Lemma 11.26]{HKT2008} it is shown how to construct elliptic curves having a stabilizer of a place of order $6$ (and hence cyclic, so the direct product of $C_2$ and $C_3$). In particular according to the proof we can see that the elliptic curve $y^2=4x^3-g_3$, $g_3 \ne 0$ satisfies this condition. It is sufficient to consider $H_1=C_2$, $H_2=C_3$. The quotient curves are clearly rational since there is ramification.

%%%%%%%%%%%%%%%%%%%%%%%%%%%%%%%%%%%%%%%%%%%%%%%%%%%
%%%%%%%%%%%%		NORM-TRACE		%%%%%%%%%%%%%%%%%%%%%%%%
%%%%%%%%%%%%%%%%%%%%%%%%%%%%%%%%%%%%%%%%%%%%%%%%%%%
\subsection{Codes from the Norm-Trace curve}\label{Ex:NormTrace}

The so called Norm-Trace curve is inspired on the surjective $\fq$-linear maps Trace 
%($Tr: \mathbb{F}_{q^{\ell}} \rightarrow\fq, x \mapsto x^{q^{\ell-1}}+x^{q^{\ell-2}}+\dots +x $) 
and Norm from $\mathbb{F}_{q^{\ell}}$ to $\fq$, 
%($N: \mathbb{F}_{q^\ell} \rightarrow \fq, x \mapsto x^{ \frac{q^\ell-1}{q-1} }$)
and generalizes the Hermitian curve obtained for $\ell=2$. This curve has been object of study in \cite{G2003} and \cite{MTT2008} for construction of codes.  The  Norm-Trace curve is defined over $\mathbb{F}_{q^{\ell}} $ by the affine equation:
$$y^{ \frac{q^\ell-1}{q-1} }=x^{q^{\ell-1}}+x^{q^{\ell-2}}+\dots +x.$$ 
 It has genus 
$g= \frac{1}{2}(q^{\ell-1}-1)(\frac{q^\ell-1}{q-1}-1)$, only one point at infinity $P_\infty$ plus $q^{2\ell-1}$ affine rational points. The functions $x$ and $y$ have pole divisor $(x)_\infty=\left(\frac{q^\ell-1}{q-1}\right)P_\infty$ and $(y)_\infty=(q^\ell-1)P_\infty$. 
Now consider two subgroups of the automorphism group of the function field $\cN$ of the curve given by
\begin{eqnarray*}
H_1&:=&\{(x, y) \mapsto (x+a, y) \mid a^{q^{\ell-1}}+a^{q^{\ell-2}}+\dots +a=0 \ \text{ and }  a \in \mathbb{F}_{q^{\ell}} \}, \text{ and }\\
H_2&:=&\{(x, y) \mapsto (x, \lambda y) \mid \lambda \in \mathbb{F}_{q^{\ell}} \text{ and } \lambda^\frac{q^\ell-1}{q-1}=1 \} ,
\end{eqnarray*}
so that  $\# \cH_1=q^{\ell-1}, \# \cH_2=\frac{q^\ell-1}{q-1}$ and the fixed fields are $\cN^{H_1}=\mathbb{F}_{q^{\ell}}(y)$, $\cN^{H_2}=\mathbb{F}_{q^{\ell}}(x)$ 
and $\cN^\cG=\mathbb{F}_{q^{\ell}}(x^{q^\ell-1}+x^{q^\ell-2}+\dots +x)=\mathbb{F}_{q^{\ell}}(y^{ \frac{q^\ell-1}{q-1} })$ by the direct product $\cG \simeq \cH_1 \times \cH_2$. The number of rational places in  $S^\cG$ that completely splits in the field extension $S| S^\cG$ is $q-1$, then we have $n=\# \cP=q^{\ell -1}(q^\ell-1)$. Following a similar construction as in section \ref{Ex:GenHermitian}
we obtain the following.

\begin{proposition}\label{normtrace}
Let $0 \leq t_i, i=1, 2$ satisfying that 
$$S=M_1q+M_2(q^\ell-1) \leq q^{\ell -1}(q^\ell-1)$$ 
for $M_1=\max\{t_1, \frac{q^\ell-1}{q-1}-2 \}$ and $M_2=\max\{t_2, q^{\ell-1} -2\}$. Choosing $m_1=\min\{t_1, \frac{q^\ell-1}{q-1}-2\}$ and $ m_2=\min\{t_2, q^{\ell-1} -2 \}$ there exists a 
$$[q^{\ell -1}(q^\ell-1), (m_1+1)(n_1+1), \geq q^{\ell -1}(q^\ell-1)-S; \frac{q^\ell-1}{q-1}-1, q^{\ell-1} -1 ]$$
recoverable code $\mathcal{C}$ over $\mathbb{F}_{q^{\ell}}$.
\end{proposition}

This codes have relative designed distance
$$\Delta(\mathcal{C})\leq \frac{S-(m_1+1)(n_1+1)+2-\left\lceil \frac{2(m_1n_1+m_1+n_1)+1}{q^\ell+q}\right\rceil}{q^{\ell -1}(q^\ell-1)}.$$

\section{General construction from subgroups with non-trivial intersection}

Here we consider non-trivial intersecting subgroups. This allows us to construct more examples.

%%%%%%%%%%%%%%%%%%%%%%%%%%%%%%%%%%%%%%%%%%%%%%%%%%%
%%%%%%%%	MAIN RESULT-NONTRIVIAL INTERSECTION	 	%%%%%%%%%%%%%%
%%%%%%%%%%%%%%%%%%%%%%%%%%%%%%%%%%%%%%%%%%%%%%%%%%%
\begin{theorem}\label{Th:GeneralConstructionNon-Trivial}
Let  $\cF | \fq$ be a function field of genus $g$.  Consider $s$ non trivial subgroups $\cH_i$ of the automorphism group of $\cF| \fq$ satisfying
\begin{enumerate}[i)]
\item[(i)] $\cG \simeq \prod_{i=1}^s\cH_i$ is a group; and let
\item[(ii)] $\#\left( \cH_i \setminus \left(\bigcup_{j < i} \cH_j\right)\right)=r_i $.
\end{enumerate}
Let $\mathcal{P}$ be a set of $n$ places in $\cF$ lying over rational places in the fixed field $\cF^{\cG}$ that are completely split in the extension $\cF | \cF^\cG$. 
Suppose that there exists a place $P_{\infty}$ of $\cF$ which is totally ramified in the extension $\cF |\cF^{\cG}$ and let $Q_{\infty}^{(i)}$ be the unique place in $\cF^{\cH_i} $ lying under $P_{\infty}$. Suppose there exists functions  $z_i, w_i$, $i=1,\ldots, s$, such that 
\begin{enumerate}[i)]
\item[(iii)] $z_i\in \cF^{\cH_i}$, $supp((z_i)_{\infty})=\{Q_{\infty}^{(i)}\}$;
\item[(iv)] $w_i\in \cF \setminus \cF^{\cH_i}$, $supp((w_i)_{\infty})=\{P_{\infty}\}$;
\item[(v)] $w_i : P^{\cH_i} \to \mathbb{F}_q$ is injective for each $i=1, \dots, s$.
\end{enumerate}
Let $t_i \geq 1$ be such that  
$$V_i := \left\{ \sum_{\ell=0}^{r_i-1} \left( \sum_{j=0}^{t_i} a_{\ell j}^{(i)}z_i^j \right)w_i^\ell \in \cF \mid  a_{\ell j}^{(i)} \in \fq \right\}\subset \mathcal{L}((n-d)P_{\infty})$$
for some $1 \leq d \leq n$ and let $V=\bigcap_{i=1}^s V_i \subset \mathcal{L}((n-d)P_{\infty})$.
If $\dim_{\mathbb{F}_q}(V)>0$, then there exists an 
$$[n,\dim_{\fq}(V),\geq d; r_1, r_2,\ldots, r_{s}]\text{-recoverable code.}$$
\end{theorem}
\proof
The proof uses  the same ideas as in Theorem \ref{Th:GeneralConstruction}.
Let $V=\bigcap_{i=1}^s V_i \subset \mathcal{L}((n-d)P_{\infty})$ and consider the following map
$$\begin{array}{llll}
e_{\mathcal{P}} \ : \ & V &\to& \mathbb{F}_q^n\\
&f &\mapsto &e_{\mathcal{P}}(f)=(f(P_1),\ldots,f(P_n)),
\end{array}$$
where $\mathcal{P}=\{P_1,\ldots, P_n\}$. The linear code $e_{\mathcal{P}}(V)$ is contained in the algebraic geometry code $C_\cL(\sum_i^n P_i, (n-d)P_{\infty})$, so it has minimum distance at least $d$. By definition the code has length $n$ and dimension $\dim_{\mathbb{F}_q}(V)$. 
Now, we deal with recoverability. Assume $P \in \mathcal{P}$ is a rational place, let $Q_i = P \cap \cF^{\cH_i}$ be the only place in $\cF^{\cH_i}$ lying under $P$, and $\cP_i \subset \cP$ be the set of places in $\cF$ over $Q_i$ for $i=1, \dots, s.$ Then $P \in \cP_i$ and $\#  \cP_i = \# \cH_i$ for $i=1, \dots , s$ by definition of the set $\cP$. 
Consider $\mathcal{Q}_i=\cP_i \setminus \left(\bigcup_{j < i} \cP_j\right)$ for any $i=1, \dots, s$. Now, $\#\mathcal{Q}_i=r_i$ and the sets $\mathcal{Q}_i$ are disjoint two by two. Using a similar approach as in Theorem \ref{Th:GeneralConstruction}, for an $f \in V_i$ one can show that $f(P)$ can be repaired by $\mathcal{Q}_i$.
\endproof

We notice that with a slight change on hypothesis of Theorem \ref{Th:GeneralConstructionNon-Trivial} we can construct codes with a different recoverability.
\begin{corollary} Whit the same notation as in Theorem \ref{Th:GeneralConstructionNon-Trivial}, changing the hypothesis $(ii)$ by 
\begin{enumerate}[i)]
\item[(ii')] there exists $m \geq 1$ sucht that $1 \leq \# \left( H_i \cap \left(\bigcup_{j \ne i} \cH_j\right)\right) \leq m,$ with $\# H_i = r_i+m$.
\end{enumerate}
we also have the existence of an
$$[n,\dim_{\fq}(V),\geq d; r_1, r_2,\ldots, r_{s}]\text{-recoverable code.}$$
\end{corollary}
\proof
The proof is similar to the one in Theorem \ref{Th:GeneralConstructionNon-Trivial}, we only need to deal with recoverability. Assume $P \in \mathcal{P}$ is a rational place, let $Q_i = P \cap \cF^{\cH_i}$ be the only place in $\cF^{\cH_i}$ lying under $P$, and $\cP_i \subset \cP$ be the set of places in $\cF$ over $Q_i$ for $i=1, \dots, s.$ Then $P \in \cP_i$ and $\#  \cP_i = \# \cH_i$ for $i=1, \dots , s$ by definition of the set $\cP$. From i) we also have $\# \left( H_i \cap \left(\bigcup_{j \ne i} \cH_j\right)\right) \leq m$ . 
For any $f \in V$, using a similar approach as in Theorem \ref{Th:GeneralConstruction} it we are going to show that  $f(P)$ can be repaired by 
$\cP_i \setminus \{P\}$ for any $i=1, \dots, s$. 
Fixing $i$ we have that $f \in V_i$ with
$\#(\cP_i \cap \left(\bigcup_{j \ne i} \cP_j\right) < m$ and $\#\cP_i=m+r_i$, so we can choose a subset 
$\mathcal{Q}_i$ of $\cP_i$ such that $\# \mathcal{Q}_i=\#\cP_i - m=r_i$ and 
$\mathcal{Q}_i \cap \left( \bigcup_{j \ne i} \cP_j \right) = \emptyset $. Then we can recover $f(P)$ from the set $\mathcal{Q}_i \subseteq \cP_i$. We also have the repair sets $\mathcal{Q}_i$ are pairwise disjoint.
\endproof

%%%%%%%%%%%%%%%%%%%%%%%%%%%%%%%%%%%%%%%%%%%%%%%%%%%
%%%%%%%%%%%%		Hermitian generalized	 			%%%%%%%%%%%%%%%%%%%%%%%%%
%%%%%%%%%%%%%%%%%%%%%%%%%%%%%%%%%%%%%%%%%%%%%%%%%%%

\subsection{Codes from the generalized Hermitian curve}
Following the same notation as in Subsection \ref{Ex:GenHermitian} we  work with the Generalized Hermitian curve $S: y^{q^\ell+1}=x^q+x$ over $\mathbb{F}_{q^{2\ell}} $. In this case, let $ \eta, \lambda \in \mathbb{F}_{q^{\ell}}^*$ with $\gcd( ord(\eta), ord(\lambda))=m>1$. Consider two subgroups of the automorphism group of the curve given by
\begin{eqnarray*}
H_1&:=&\{\sigma_i\,:\, (x, y) \mapsto (x+a, \eta^i y) \mid a^q+a=0 \ \text{ and }  a \in \mathbb{F}_{q^{2\ell}}, i=0, \dots, ord(\eta)-1 \}, \text{ and }\\
H_2&:=&\{\sigma_i\,:\,(x, y) \mapsto (x, \lambda^i y) \mid \lambda \in \mathbb{F}_{q^{2\ell}}, i=0, \dots, ord(\lambda)-1 \} ,
\end{eqnarray*}
so that  $\# \cH_1=q \cdot ord(\eta), \# \cH_2=ord(\lambda)$, they comute and $\cG = \cH_1 \cH_2$ has order $\frac{q \cdot ord(\eta)ord(\lambda) }{m}$.
The fixed fields are 
$$ 
\cS^{H_1}=\mathbb{F}_{q^{2\ell}}(x^q+x, y^{\frac{q^\ell+1}{ord(\eta)}}), \, \cS^{H_2}=\mathbb{F}_{q^{2\ell}}(x, y^{\frac{q^\ell+1}{ord(\lambda)}}), \text{ and } \cS^\cG=\mathbb{F}_{q^{2\ell}}(x^q+x, y^{\frac{q^\ell+1}{m}})
$$ 
wher $\cG \simeq \cH_1 \times \cH_2$ is the direct product. The number of rational places in  $S^\cG$ that completely splits in the field extension $S| S^\cG$ is $m(q^\ell-1)$, then we have $n=\# \cP=(q^{\ell}-1)ord(\eta)ord(\lambda)$.
Let  $ 0 \leq t_i, i=1, 2$ and consider the vector spaces
\begin{eqnarray*}
	V_1 &:=&\left\{\sum_{\ell=0}^{ord(\eta)-2} \left(\sum_{j=0}^{t_1 }a_{\ell, j}y^{j\frac{q^\ell+1}{ord(\eta)} }\right)x^\ell   \ : \ a_{\ell, j} \in \mathbb{F}_{q^{2\ell}}\right \},\\
	V_2 &:=&\left\{\sum_{\ell=0}^{ord(\lambda)-2} \left(\sum_{j=0}^{t_2 }a_{\ell, j}x^j  \right)y^\ell   \ : \ a_{\ell, j} \in \mathbb{F}_{q^{2\ell}}\right \},\\
\end{eqnarray*}
Choose $t_i$ such that for $M_1=\max\{t_1, ord(\lambda)-2 \}, M_2=\max\{t_2, ord(\eta)-2 \}$  we have $$S=\max\{M_1q\frac{q^\ell+1}{ord(\eta)}+(q^\ell+1)(ord(\eta)-2), M_2(q^\ell+1)+q(ord(\lambda)-2) \} \leq n.$$
 Then  $V_i \subset \cL(SP_\infty),$ for $ i=1, 2$, and for  $m_1=\min\{t_1, ord(\lambda)-2 \}, m_2=\min\{t_2, ord(\eta)-2 \}$ we obtain
$$V=V_1\cap V_2=\left\{\sum_{\ell=0}^{m_1} \left(\sum_{j=0}^{m_2 }a_{\ell, j}x^j  \right)y^\ell   \ : \ a_{\ell, j} \in \mathbb{F}_{q^{2\ell}}\right \}  \textrm{ and } \dim _{\mathbb{F}_{q^2}}V=(m_1+1)(m_2+1).$$
The minimum distance of the code satisfies $d \geq q^{2\ell +1}-q- S$. We have obtained the following proposition.

\begin{proposition}\label{genherm1}
Let $q$ be odd, $ \eta, \lambda \in \mathbb{F}_{q^{\ell}}^*$ with $\gcd( ord(\eta), ord(\lambda))=m>1$ and $n=(q^{\ell}-1)ord(\eta)ord(\lambda)$. 
Choose $t_i$ such that for $M_1=\max\{t_1, ord(\lambda)-2 \}, M_2=\max\{t_2, ord(\eta)-2 \}$  we have $$S=\max\{M_1q\frac{q^\ell+1}{ord(\eta)}+(q^\ell+1)(ord(\eta)-2), M_2(q^\ell+1)+q(ord(\lambda)-2) \} \leq n.$$
$m_1=\min\{t_1, ord(\lambda)-2 \}, m_2=\min\{t_2, ord(\eta)-2 \}$ there exists a 
$$[(q^{\ell}-1)ord(\eta)ord(\lambda), (m_1+1)(n_1+1), \geq n-S; q\cdot ord(\eta), ord(\lambda)-m]$$
recoverable code $\mathcal{C}$ over $\mathbb{F}_{q^{2\ell}}$.
\end{proposition}

The codes arising from Proposition \ref{genherm1} have relative designed distance
$$\Delta(\mathcal{C})\leq \frac{S-(m_1+1)(n_1+1)+2-\left\lceil \frac{2(m_1n_1+m_1+n_1)+1}{q\cdot ord(\eta)+ord(\lambda)-m+1}\right\rceil}{(q^{\ell}-1)ord(\eta)ord(\lambda)}.$$

\begin{remark}
Every example proposed in the previous section can be used to apply the construction introduced in this section as well obtaining different parameters. Indeed it is sufficient to enlarge a bit one of the two groups including elements from the second one.
\end{remark}

\section*{Acknowledgement}
This work was partially done during the visit of the third author to the University of Perugia as Visiting  Researcher (project name: ``Algebraic Curves and Applications", 23/01/2019--23/02/2019) and the visit of the first author to the Federal University of Rio de Janeiro(Capes - Proex).

\end{document}